\documentstyle[12pt]{article}
\newcommand{\ul}{\underline}
\begin{document}

\title{Topological Hilbert Nullstellensatz for Bergman Spaces}
\author{R{\u{a}}zvan Gelca}
\maketitle

{\centerline {\bf 1. INTRODUCTION}}

\bigskip

In the present paper we intend to discuss problems related to the
classification of invariant subspaces for multiplication
operators on spaces of functions of several variables. Unlike
the one-dimensional situation, the several variables case proves to be
very hard, and little progress has been made until now.
In [4] and [6] the attention was focused on invariant
subspaces that arise as closures of ideals from some dense ring,
the main result proved in this direction being the rigidity theorem
(see [4]). In order to understand this kind of invariant subspaces,
one needs first to understand their corresponding ideals. In particular
one needs to classify the ideals that are closed in the relative topology
induced on the ring by the space of functions.

In [3] it was conjectured that an ideal of polynomials is
closed in the relative topology induced by the Hardy space of the unit
polydisk if and only if every irreducible component of the variety
of the ideal intersects the polydisk. In [5] the conjecture was
 proved in dimension 2. The idea of the proof is to relate the
conjecture to a topological version of the Hilbert Nullstellensatz for
the ring of polynomials. In this paper we intend to prove similar
results for Bergman space topologies. One should note that our methods
provide new proofs for the one-dimensional case as well. Some
of the techniques come from commutative algebra, and our basic source is
 [1].

Let us start by recalling some facts from [5].
 Let ${\cal R}$ be a commutative Noetherian ring with unit and let $\tau $
be a topology on ${\cal R}$ for which the addition is continuous
and multiplication is separately continuous in each variable. Examples of
such rings are the ring of polynomials with the topology induced by
the Hardy space or by some Bergman space, the ring
 of germs of analytic functions in a neighborhood of the closure
of some domain satisfying certain properties [7] with the topology
induced by the Bergman space, etc.

In [5] it has been shown that the radical of a closed ideal is
closed, and that the prime ideals associated to a closed ideal are
also closed. Define ${\cal C}$ to be the set of all maximal ideals
${\cal M}$ with the property that ${\cal M}^n$ is closed for all $n\in
{\bf N}$.
In our examples ${\cal C}$ will consist of the maximal ideals that
correspond to points in the domain on which the spaces of functions are
defined. The pair $({\cal R}, \tau)$
is said to satisfy the topological Hilbert Nullstellensatz if every
ideal $I\subset {\cal R}$ is either dense or there exists ${\cal M}
\in {\cal C}$ with $I\subset {\cal M}$. The Nullstellensatz is easy to
prove in the one-dimensional case, but becomes hard when passing to
several variables. In what follows we show that it holds for various
topologies in the two variable case. We will need the following
elementary result, whose proof can be found in [5].

\medskip

\underline{Proposition 1.1.} Let $p(z)=a(z-z_1)(z-z_2)\cdots (z-z_m)$ be
a polynomial with the property that $|z_i|\geq 1$, $i=1,2,\cdots, m$.
Then for all $1/2<r<1$ and $z$ with $|z|\leq 1$ we have $|p(z)/
p(rz)|\leq 2^m$.

\bigskip

{\centerline {\bf 2. THE RING OF POLYNOMIALS}}

\bigskip

Let $\Omega\in {\bf C}^n$ be a bounded complete Reinhardt domain (i.e. a
domain with the property that for all $x\in \Omega$ and $|z|\leq 1$,
$zx$ is in $\Omega$), and let
$L_a^2(\Omega )$ be the Bergman space of $\Omega $. For a polynomial
$p \in {\bf C}[z_1,z_2,\cdots ,z_n]$ we denote by $d(p)$ the sum of
the degrees of $p$ in each variable.

\medskip

\underline{Lemma 2.1.} Let $p \in {\bf C}[z_1,z_2,\cdots ,z_n]$ be a
polynomial having no zeros in $\Omega$. Then for every $1/2<r<1$ and
$(z_1,z_2,\cdots ,z_n)\in \Omega $, $|p(z_1,z_2,\cdots ,z_n)/
p(rz_1,rz_2,\cdots ,rz_n)|\leq 2^{d(p)}$.

\medskip

\underline{Proof:}  The proof will be done by induction on the number of
variables. By Proposition 1.1 the property is true for $n=1$.

Fix $(z_1,z_2,\cdots ,z_n)\in \Omega $. Let $D=\{ w, \
(w,z_2,\cdots ,z_n)\in \Omega\}$. Then $D$ is a disk and
$p(w,z_2,\cdots ,z_n)$
has no zeros in $D$. Thus $|p(z_1,z_2,\cdots ,z_n)/
p(rz_1,z_2,\cdots ,z_n)|\leq 2^{d_1}$, where we denote by $d_i$ the degree
of $p$ in the variable $z_i$. Applying the induction hypothesis
to $p(z_1,\cdot, \cdots ,\cdot)$ we get $|p(rz_1,z_2,\cdots ,z_n)/
p(rz_1,rz_2,\cdots ,rz_n)|\leq 2^{d_2+\cdots d_n}$, thus
$|p(z_1,z_2,\cdots ,z_n)/
p(rz_1,rz_2,\cdots ,rz_n)|\leq 2^{d(p)}$.$\Box$

\medskip

\ul{Proposition 2.2.} Let $p \in {\bf C}[z_1,z_2,\cdots ,z_n]$ be
a polynomial having no zeros in $\Omega$. Then $pL^2_a (\Omega )$ is dense
in $L^2_a (\Omega )$.

\medskip

\ul{Proof:} Let $f_r(z_1,z_2,\cdots ,z_n)=p(z_1,z_2,\cdots ,z_n)/
p(rz_1,rz_2,\cdots ,rz_n)$. Then for $1/2<r<1$, $f_r \in pL^2_a (\Omega )$.
By Lemma 2.1 the family $f_r$ is bounded, and since for $r\rightarrow 1$, $f_r$
converges uniformly on compact subsets of $\Omega$
 to $1$, it follows that it converges in
$L^2_a (\Omega )$ to $1$. The conclusion follows by using the fact that
the ring of polynomials is dense in $L^2_a (\Omega )$.$\Box$

\medskip

\ul{Corollary 2.3.} The topological Hilbert Nullstellensatz holds for
principal ideals in the ring of polynomials with the topology induced by
$L^2_a(\Omega )$.

\medskip

\ul{Proof:} If ${\cal M}$ is a maximal ideal that corresponds
to a point in $\Omega $ then ${\cal M}^n$ is closed for all $n$,
thus ${\cal M}\in {\cal C}$. This shows that if
the variety of an ideal intersects $\Omega $ then the ideal is contained
in some maximal ideal from ${\cal C}$. If the variety of a principal
ideal does not intersect $\Omega $, then the ideal is dense and the conclusion
follows.$\Box$

\medskip

We now exhibit some examples of topologies for which the topological
Hilbert Nullstellensatz will hold for the ring of polynomials in two
variables.
Following [2] we define

\hspace{5mm} $\Omega _{p,q}:=\{(z_1,z_2)\in {\bf C}^2 \ \ |
 \ \ |z_1|^p+|z_2|^q<1\}.$

where $0<p,q<\infty $. Note that $\Omega _{2,2}$ is the unit ball in ${\bf C}^2
$.

\medskip

\ul{Theorem 2.4} The ring ${\bf C}[z_1,z_2]$ with the topology induced by
$L^2(\Omega _{p,q})$ satisfies the topological Hilbert Nullstellensatz.

\medskip

\ul{Proof:} By Proposition 4.1 in [5] ${\cal C}=\Omega$ and the
maximal ideals not contained in ${\cal C}$ are dense, thus
 the Nullstellensatz holds for
maximal ideals. The previous corollary shows that it also
holds for principal ideals. Since the ring ${\bf C}[z_1,z_2]$ has
dimension 2, a prime ideal is either maximal or principal, hence the
Nullstellensatz holds for prime ideals.

Now let $I$ be an arbitrary ideal. Then  $I$ is either
contained in an element from $ {\cal C}$, or its variety does not
intersect $\Omega $. In the latter case all prime ideals associated to
$I$ are dense. Let $P_1,P_2, \cdots ,P_m$ be these ideals.
By Noetherianess (see Proposition 7.14 in [1]) there is a $k$ such that
$(P_1P_2\cdots P_m)^k\subset I$, and since $(P_1P_2\cdots P_m)^k$ is dense
it follows that $I$ is dense. This proves the theorem.$\Box$

\medskip

By using Theorem 1.3 in [5] we get the following analogue of the
Douglas-Paulsen conjecture.

\medskip

\ul{Corollary 2.5.} Let ${\bf C}[z_1,z_2]$ be endowed with the topology
induced by $L^2_a(\Omega _{p,q})$. Then an ideal is closed if and only
if each of the irreducible components of its zero set intersects
$\Omega _{p,q}$.

\bigskip

{\begin{center}
{\bf 3. THE RING OF ANALYTIC FUNCTIONS IN A NEIGHBORHOOD
OF THE UNIT BALL}
\end{center}}

\bigskip

Let us denote by ${\bf B}^2$ the open unit ball in ${\bf C}^2$ and
by ${\cal O}(\overline{{\bf B}^2})$ the ring of germs of analytic functions
defined in a neighborhood of $\overline{{\bf B}^2}$. We want to prove the
topological Hilbert Nullstellensatz for ${\cal O}(\overline{{\bf B}^2})$
with the topology induced by the Bergman space. For $k\geq 1$,
denote  by $S^k$
the unit sphere in ${\bf R}^{k+1}$ and by ${\bf D}$ the unit disk in
the plane.

\medskip

\ul{Lemma 3.1.} Let $f\in {\cal O}(\overline{{\bf B}^2})$, $f$ not identically
equal to zero. Then there exists a transformation $ \rho :{\bf C}^2
\rightarrow {\bf C}^2$ such that $f \circ \rho (0,e^{i \alpha})\neq 0$
for all $\alpha \in [0,2\pi )$.

\medskip

\ul{Proof:} For $S^3=\partial {\bf B}^2$ consider  the Hopf fibration
\begin{eqnarray}
S^1\rightarrow S^3 \stackrel{\pi}{\rightarrow} S^2
\end{eqnarray}
where we recall that the projection $\pi $ is given by the
equivalence relation $(a,b)\sim (\lambda a,\lambda b)$ for $|\lambda |=1$.
If we denote by $V(f)$ the zero set of $f$, then $V(f)\cap S^3$ has
dimension at most one, hence $\pi (V(f)\cap S^3)$ has dimension
at most one. It follows that there exists $x\in S^2\backslash
\pi (V(f)\cap S^3)$. Let $(a,b)\in {\bf C}^2$, $|a|^2+|b|^2=1$ with the
property that $\pi (a,b)=x$. Then $f(ae^{i\alpha },be^{i\alpha })\neq
0$ for all $\alpha \in [0,2\pi )$. If we choose $\rho =
\left(
\begin{array}{clcr}
\bar{b} & a \\
-\bar{a} & b
\end{array}
\right)$
then $f\circ \rho$ satisfies the desired property.$\Box$

\medskip

\ul{Proposition 3.2.} Let $f\in {\cal O}(\overline{{\bf B}^2})$
be an analytic function having no zeros in ${\bf B}^2$. Then the
space $fL^2_a({\bf B}^2)$ is dense in $L^2_a({\bf B}^2)$.

\medskip

\ul{Proof:} By Lemma 3.1 we may assume that if $ (w_1,w_2)\in
\overline{{\bf B}^2}$
and $w_1=0$ then $f(w_1,w_2)\neq 0$. Let us show that there
exists $C>0$ such that for every $1/2<r<1$ and $(z_1,z_2)\in
\overline{{\bf B}^2}$,
$|f(z_1,z_2)/f(rz_1,z_2)|<C$.

Fix $(w_1,w_2)\in \overline{{\bf B}^2}$, and assume that $w_1\neq 0$. Note
that $f(z,w_2)$ is not identically zero as a function of $z_1$,
otherwise $(0,w_2)$ would be a zero for $f$ lying inside of ${\bf B}^2$.
Thus there exists $a>0$ such that $w_1\in a{\bf D}$ and $f(\cdot ,w_2)$
is defined in a neighborhood of $\overline{a{\bf D}}$
 and has no zeros on $aS^1=
\{z\ | \ |z|=a\}$. It follows that $f(\cdot ,w_2)$ has a finite number of
zeros in $a{\bf D}$, say $m$, counting multiplicities. By Rouch{\'{e}}'s
Theorem there is a compact neighborhood $K$ of $w_2$ such that
$\overline{a{\bf D}}\times K$ is contained in the domain of $f$, and for
every $z_2\in K$, $f(\cdot ,z_2)$ has exactly $m$ zeros in $a{\bf D}$,
counting multiplicities, and no zero on $aS^1$.

Thus on $\overline{a{\bf D}}\times K$ we can write $f(z_1,z_2)=
p_{z_2}(z_1)g_{z_2}(z_1)$ where for each $z_2$, $p_{z_2}(z_1)$ is a
polynomial of degree $m$ and $g_{z_2}(z_1)$ is an analytic function
having no zeros in $\overline{a{\bf D}}$. Another application of
Rouch{\'{e}}'s Theorem and the maximum modulus principle shows that
$g_{z_2}$ depends continuously on $z_2$.

It follows that for $1/2\leq r \leq 1$ the family $\{g_{z_2}(rz_1)\}_r$
is bounded away from zero, hence
\begin{eqnarray}
C_1=sup_{1/2\leq r\leq 1}sup_{\overline{a{\bf D}}\times K}|g_{z_2}(z_1)/
g_{z_2}(rz_1)|<\infty .
\end{eqnarray}

By Proposition 1.1, for $1/2\leq r \leq 1$
and $(z_1,z_2)\in \overline{a{\bf D}}\times K
\cap
\overline{{\bf B}^2}$
\begin{eqnarray}
|p_{z_2}(z_1)/p_{z_2}(rz_1)|\leq 2^m.
\end{eqnarray}

Thus there exists a neighborhood $U$ of
$(w_1,w_2)$ and a constant $C_2>0$ such that for $1/2<r<1$ and $(z_1,z_2)\in U
\cap
\overline{{\bf B}^2}$
\begin{eqnarray}
|f(z_1,z_2)/f(rz_1,z_2)|<C_2.
\end{eqnarray}

If $w_1=0$ then $f(z_1,z_2)\neq 0$ in a neighborhood of $(w_1,w_2)$, thus
a similar inequality holds there. From the compactness of $\overline{{\bf
 B}^2}$
it follows that there exists a constant $C>0$ such that for $1/2<r<1$
and $(z_1,z_2)\in \overline{{\bf B}^2}$, $|f(z_1,z_2)/f(rz_1,z_2)|<C$.

As in the proof of Proposition 2.2, the family $h_r(z_1,z_2)=
f(z_1,z_2)/f(rz_1,z_2)$ is in $fL^2_a({\bf B}^2)$ and tends to
$1$ as $r\rightarrow 1$, so the conclusion follows.$\Box$

\medskip

\ul{Theorem 3.3.} The ring ${\cal O}(\overline{{\bf B}^2})$ with the topology
induced by $L^2_a({\bf B}^2)$ satisfies the topological Hilbert
Nullstellensatz.

\medskip

\ul{Proof:} The ring ${\cal O}(\overline{{\bf B}^2})$ is Noetherian [7],
and has dimension 2. Indeed, if there existed distinct prime
ideals $P_0\subset P_1\subset P_2\subset P_3$, by localizing at a
maximal ideal ${\cal M}\supset P_3$ we would get a chain of four
distinct prime ideals in the local ring ${\cal O}_{\cal M}$, which
would contradict the fact that the latter ring has dimension 2.
So the proof of Theorem 2.4 applies mutatis mutandis to give
the desired conclusion.$\Box$

\medskip

\ul{Corollary 3.4.} Let ${\cal O}(\overline{{\bf B}^2})$ be endowed
with the topology induced by the Bergman space. Then an ideal is
closed if and only if each irreducible component of its zero set
intersects the unit ball.

\bigskip

{\bf References}

\medskip

\noindent [1] M. F. ATIYAH, I. G.  MACDONALD, Introduction to Commutative
 Algebra, Addison-Wesley, 1969.

\noindent [2] R. CURTO, N. SALINAS, Spectral properties of cyclic subnormal
m-tuples, {\em Amer.
J. Math.}, {\bf 1}(1985), 113--138.

\noindent [3] R. G. DOUGLAS, V. PAULSEN, Hilbert Modules over Function
Algebras, Longman
Scientific and Technical, London, 1989.

\noindent [4] R. G. DOUGLAS, V. PAULSEN, C. SAH, K. YAN,  Algebraic reduction
and rigidity
for Hilbert modules, {\em Amer. J. Math.}, {\bf 117}(1995),
 73--93.

\noindent [5] R. GELCA, Rings with topologies induced by spaces of functions,
{\em Houston J. Math},
{\bf 2}(1995), 395--405.

\noindent [6] V. PAULSEN,  Rigidity theorems in spaces of analytic functions,
{\em Proc.
Symp. Pure
Math.}, {\bf  51}(1991), Part 1.

\noindent [7] Y. T. SIU, Noetherianess of rings of holomorphic functions
on Stein compact
 sets, {\em Proc. Amer.
 Math. Soc.}, {\bf 21}(1969), 483--489.

\medskip

\noindent Department of Mathematics,
The University of Iowa, Iowa City, IA 52242

\noindent (mailing address)

\noindent {\em E-mail: rgelca@math.uiowa.edu}

\medskip

and

\medskip

\noindent Institute of
Mathematics of the Romanian Academy, P.O.Box 1-764, Bucharest
70700, Romania.

\end{document}